\documentclass{article}
\usepackage{amsmath,amssymb}
\usepackage{amsthm}
\usepackage{amsfonts}
\usepackage{graphicx}

\newcommand{\tmtextbf}[1]{{\bfseries{#1}}}

\newcommand{\eps}{{\epsilon}}
\newcommand{\var}{{\hbox {var}}}

\newtheorem{theorem}{Theorem}

\newtheorem {prop}[theorem]    {Proposition}

\begin{document}

\title{Asymmetry of near-critical percolation interfaces}
\author{Pierre Nolin \and Wendelin Werner}
\date {Universit\'e Paris-Sud 11 and \'Ecole Normale Sup\'erieure\footnote {Research supported by the Agence Nationale pour la Recherche under the grant ANR-06-BLAN-0058.}}
\maketitle

\begin{abstract}
  We study the possible scaling limits of percolation interfaces in two
  dimensions on the triangular lattice. When one lets the percolation
  parameter $p(N)$ vary with the size $N$ of the box that one is
  considering, three possibilities arise in the large-scale limit. 
  It is known that when $p(N)$ does not converge to $1/2$ fast enough, then the scaling limits are degenerate, whereas if $p(N) - 1 / 2$ goes to zero quickly, the scaling limits are SLE(6) as when $p=1/2$. We study some properties of the (non-void) intermediate regime where the large scale behavior is neither SLE(6) nor degenerate. 
We prove that in this case, the law of any scaling limit is singular with respect to that of
  SLE(6), even if it is still supported on the set of curves with Hausdorff dimension equal to $7/4$.
\end{abstract}

\section{Introduction}

In this paper, we study site percolation on the triangular planar lattice.
Recall that this can be viewed as a random coloring of the hexagonal cells of
a honeycomb lattice, where the color (black or white) of each cell is chosen
independently of the others: each of these cells has a probability $p$ to be black and $1 -
p$ to be white, for some parameter $p$ between $0$ and $1$. In percolation
theory, one is interested in the connectivity properties of the set of black hexagons (or the set of white ones). They can be regrouped into connected components (or clusters). The
phase transition for percolation on this lattice occurs at $p = 1 / 2$. Often,
it is described mathematically as follows, in terms of almost sure properties
of percolation in the infinite lattice: when $p < 1 / 2$, there exists with
probability $1$ no infinite cluster of black sites (subcritical regime) and an
infinite cluster of white sites, and conversely when $p > 1 / 2$, there is an
infinite cluster of black sites (supercritical regime) but no infinite cluster
of white sites. In the critical case where $p = 1 / 2$, there exists neither
an infinite white cluster, nor an infinite black cluster -- but if one takes a
finite large piece  $\Lambda$ of the lattice, one will see white and black clusters of
size comparable to that of $\Lambda$. See e.g. {\cite{Grbook,Kebook}} for an
introduction to percolation.

A lot of progress has been made recently in the understanding of the
large-scale behavior of critical percolation: in particular, Smirnov {\cite{Sm}} proved  conformal invariance of
the connection probabilities, which allowed to make the link
{\cite{Sm,CN1,Sm2}} with the Schramm-Loewner Evolution (SLE) with parameter
6 introduced in {\cite{Sch}}, and to use the SLE technology and computations
{\cite{LSW1,LSW2}} to derive further properties of critical percolation, such
as the value of some critical exponents, describing the asymptotic behavior
of the probabilities of certain exceptional events (arm exponents)
{\cite{LSW5,SmW}}. We refer to {\cite{Wpc}} for a survey.

\begin{figure}
\begin{center}
\includegraphics[width=9cm]{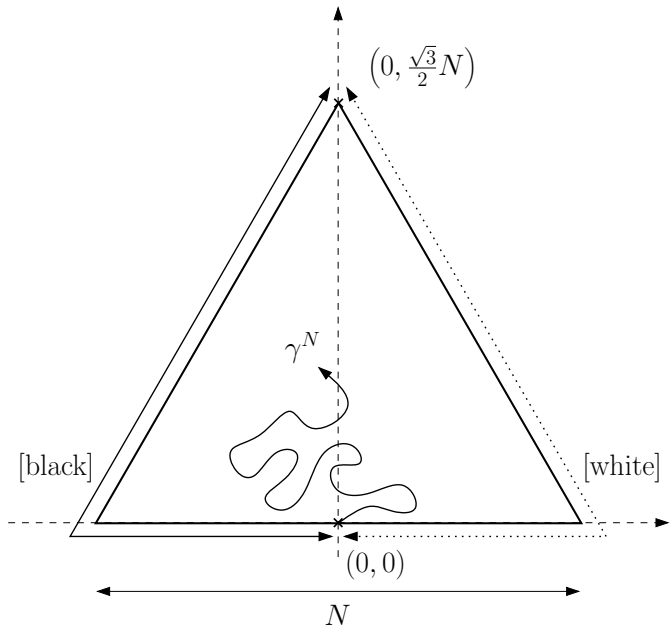}
\caption{The triangle $T_N$ and the interface $\gamma^N$ (sketch).}
\label{triangle}
\end{center}
\end{figure}

One precise relation to SLE goes as follows: we consider the large equilateral
triangle $T_N$ with even side length $N$ on the triangular grid such that the
middle of the bottom part is the origin and the top point is the point at
distance $\sqrt{3} N / 2$ above the origin. We decide to color all cells on
the boundary of the triangle, in white if their $x$-coordinate is positive and
in black if their $x$-coordinate is negative, and we perform critical
percolation in the inside of $T_N$. Then, we consider the interface $\gamma^N$ (viewed as a path on the hexagonal lattice dual to the triangular lattice) between the set of black sites attached to the left part of the triangle
and the set of white sites connected to the right part of the triangle (see Figure \ref{triangle}).
When $N \to \infty$, the law of the rescaled interface $\Gamma_N:= \gamma^N / N$ converges (in an appropriate
topology) to that of the SLE(6) process from $(0, 0)$ to $(0, \sqrt{3} / 2)$
in the equilateral triangle with unit side length. See {\cite{CN2,Sm2,Wpc}} for
details, and e.g. {\cite{Lbook}} for an introduction to SLE. Thanks to this
convergence result, one is able to deduce properties of critical percolation
from the properties of SLE. For instance (and we shall come back to this
later), one can prove that the typical number of steps of the path $\gamma^N$
is of the order $N^{7 / 4}$ (more precisely, for each $\epsilon > 0$, the
probability that the number of steps is between $N^{7 / 4 - \epsilon}$ and
$N^{7 / 4 + \epsilon}$ goes to $1$ as $N \to \infty$).

Understanding the behavior of critical percolation allows also to derive some
properties of percolation when the parameter $p$ is very close to $1 / 2$.
These are the scaling (and hyperscaling) relations that were first
developed in the physics literature (the names of Fischer, Kadanoff and Wilson are often associated to these ideas),
 and later rigorously derived in the case
of percolation by Kesten {\cite{Ke4}} (see also {\cite{N2,Wpc}}). Before
explaining these relations in a little more detail, let us first make a
trivial observation to help the newcomer to get a feeling of what goes on.
Suppose that we consider the percolation interface in the large equilateral
triangle as before, but choose $p = p (N)$ in such a way that $|p - 1 / 2| = o
(N^{-7 / 4 - \epsilon})$. Then, one can couple it with a critical
percolation interface (i.e. with the case when $p = 1 / 2$ exactly) in the
same triangle in such a way that the two paths are identical except on an
event of probability that goes to $0$ as $N\to \infty$ (this is just because
the probability that one of the neighbors of the interface changed its color
is bounded by $(p - 1 / 2)$ times the number of neighbors of the interface --
and we have just seen that this number is $o (N^{7 / 4 + \epsilon})$).
Hence, if one lets $N$ go to infinity, the scaling limit of the interface is
still SLE(6). In particular, the probability that the interface $\Gamma_N$
hits the right side of the triangle before the left side of the triangle -- let us call $R (p, N)$ this probability -- goes to $1 / 2$ as $N \to \infty$.

On the other hand, for any fixed $\epsilon \in (0,1/2)$, one can define
\[ p^{\ast} (N) = p^{\ast} (N,\epsilon) := \inf \{p : R (p, N) > 1 / 2 +
   \epsilon\}. \]
For this choice of $p = p^{\ast} (N)$, if one looks at the possible limiting
behavior of $\Gamma_N$, it is clear that the law can not be exactly SLE(6)
anymore because it will hit the right side of the triangle before the left one 
with probability at least  $1 / 2+\eps$. It is therefore natural to ask what can happen
to the scaling limit of this curve when $N \to \infty$ 
in this regime, and to see how it is
related (or not) to SLE(6).

One can equivalently define the so-called correlation length $L (p) = L (p, \epsilon)$ in
such a way that $p^{\ast} (L(p),\epsilon) \simeq p$. In other words, for $p
> 1 / 2$,
\[ L (p) = L (p, \epsilon) := \inf \{N : R (p, N) > 1 / 2 + \epsilon\} \]
(note that for $p>1/2$ fixed, $R (p, N) \to 1$ as $N \to \infty$). Kesten {\cite{Ke4}} has shown that it is possible to
deduce from the arm exponents of critical percolation the behavior of $L (p)$ as
$p \to 1 / 2$. This derivation relies on the ``four-arm exponent'': the crossing probability increases exactly when ``pivotal'' sites are flipped. Combining Kesten's results with the
exponents computed using SLE, one gets {\cite{SmW}} that
\[ L (p) = (p - 1 / 2)^{- 4 / 3 + o (1)} \]
when $p \to 1/2^+$, for any fixed choice of $\epsilon \in (0,1/2)$.
This is (see \cite {Ke4, N2, Wpc}) a crucial step in the rigorous proof of the fact that 
the ``density'' $\theta (p)$ of the infinite cluster for percolation with parameter $p$ decays like 
$(p-1/2)^{5/36 + o (1)}$
as $p$ decays to $1/2^+$.

We shall see (and this is quite easy) that in order to get a non-trivial limit for $\Gamma_N$ (i.e. neither
SLE(6) nor a path that just sticks to the boundary of $T$), one has to take $p (N)$ in such a way that $N$ is of the
order of the correlation length $L (p)$ i.e. that $p (N) \in [p^{\ast}
(N,\epsilon), p^{\ast} (N,\epsilon')]$ for some $\epsilon < \epsilon'$. The terminology to describe this regime depends on the circumstances and the papers:
near-critical, off-critical,  finite-size scaling. It is also very closely related to the scaling and hyperscaling relations.  Anyway, it has been the subject of numerous and interesting works, see e.g. {\cite{A, Ch, BCKS2,CCFS}} and the references therein.

One general rule behind the derivation of these results is that as long as $n \le L(p)$, the behavior remains roughly the same as at criticality. We will remind some related results in the next sections, but for example, the probabilities of existence of crossings of annuli for near-critical percolation are bounded by constant times those at criticality. This implies that the  exponents describing critical percolation also describe off-critical percolation and lead to Kesten's scaling relations. In \cite {CFN1, CFN2}, Camia, Fontes and Newman have suggested that the scaling limit of the near-critical picture could be obtained from that of the critical picture by just opening a Poissonian family of ``pivotal'' points. We understand that \cite {GPS2} have succeeded in proving this fact. A consequence would be that the subsequential limits that we are discussing in the present paper are in fact limits, and that they are all related to each other by scaling. 
Let us also mention that near-critical percolation is closely related to the questions concerning the noise-sensitivity of percolation, such as studied in \cite {SS, GPS1}.  

Before describing the results of the present paper, let us just emphasize the relevance of such near-critical models for physical applications.
Suppose for instance that one is considering a percolation model in a large box, and that one ``heats'' progressively the system i.e. that cells become progressively black one after the other (and independently). Then, we know that the percolation parameter at which a left-to-right crossing of the big box will appear is (with a large probability) close to the critical one. But, if one looks at it more precisely, it is easy to see that the value at which the crossing occurs will be $p^* (N,u)$ where $u$ is a uniform random variable in $[-1/2,1/2]$. Hence, at this moment, the picture and its fine properties are that of near-critical percolation. Another example where this near-critical picture appears spontaneously is that of gradient percolation, an inhomogeneous percolation model introduced in \cite {SG} and studied mathematically in \cite {N1}.

\begin{figure}
\begin{center}
\includegraphics[width=4in]{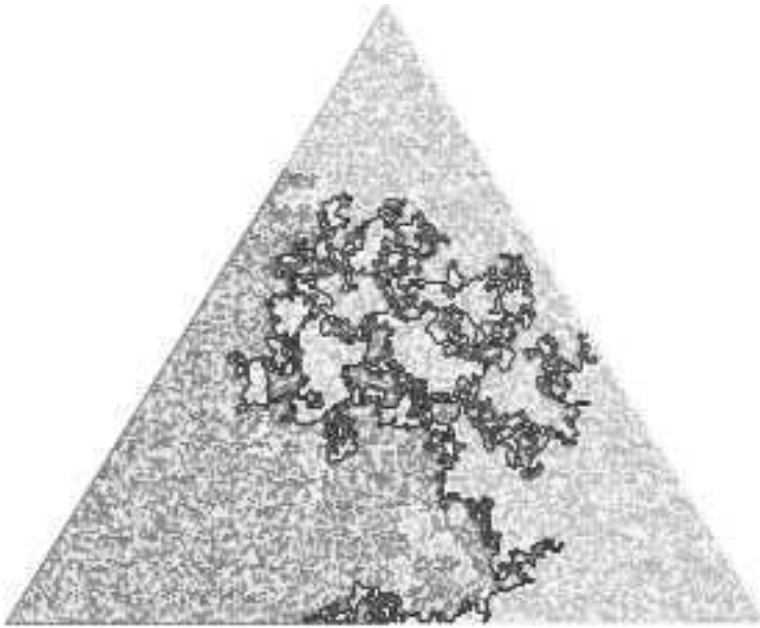}
\caption{A coupling of a near-critical interface with critical percolation.}
\label{nicepic}
\end{center}
\end{figure}

Our focus here will be on the scaling limits of near-critical interfaces. 
Here is a list of results that we shall derive in the present paper.
Choose a sequence $(p(N), N \ge 1)$ and study the behavior of the law of $\Gamma_N$
as $N \to \infty$:
\begin {itemize}
\item
For any $\epsilon \in (0,1/2)$, if we choose $p(N)$ in such a way that 
$p(N)$ is close to $p^* (N, \eps)$, then there exist subsequential limits for the law of $\Gamma_N$.
\item
The limiting laws are all singular with respect to the law of SLE(6).
\item 
 The only other possible scaling limits are SLE(6) itself (when $p(N)$ is too close to $1/2$) or 
a degenerate case where $\gamma$ follows (the right part of) the boundary of the triangle (this is when $p(N)$ is too large).
\item
The critical exponents associated to these non-degenerate scaling limits are the same as those of SLE(6). In particular, the Hausdorff dimension 
of the curves is almost surely $7/4$.
\end {itemize}

\begin{figure}
\begin{center}
\includegraphics[width=4in]{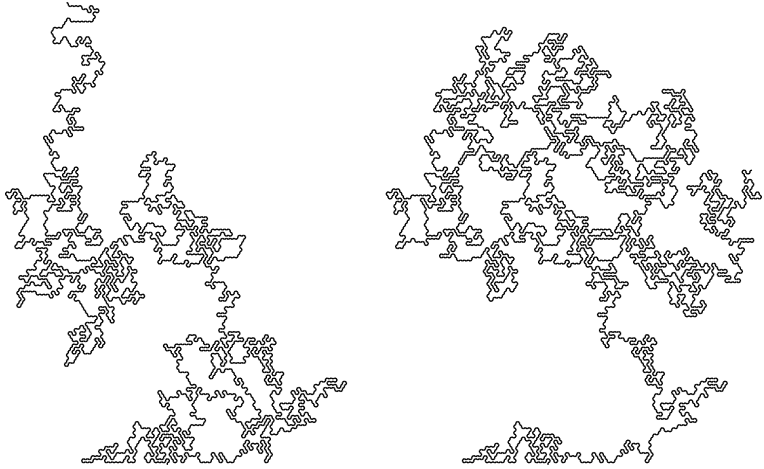}
\caption{The interfaces at $p_c$ and at $p^* (N, \epsilon)$ of Figure \ref {nicepic} depicted separately.}
\label{2interfaces}
\end{center}
\end{figure}

The first result is a rather direct application of the tightness arguments of Aizenman and Burchard \cite {AB} and the definition of the correlation length.
The last two statements follow from the ideas of Kesten's paper \cite {Ke4}.
 Our main result is probably the second one, the fact that the laws of non-trivial scaling limits are singular with respect to SLE(6).
One way to explain it is the following. In the off-critical
regime, one sees on a macroscopic scale a difference between the law of the
interface and that of the critical percolation interface (i.e. the
non-critical interface is more to the ``right'' for instance). If one zooms in
by a factor $\lambda$, one still sees a difference, but this difference tends
to disappear, because one is not looking at a picture of size ca. $L (p)$ any
more, but at a picture of size $L (p) / \lambda$. The question is whether this
difference disappears sufficiently fast when $\lambda \to \infty$ or not. Note
that one can show that (just as for the critical interface) the number of
boxes of size $N / \lambda$ visited by the path is of order $\lambda^{7 / 4}$
when $\lambda$ is large (and $N$ very large). Either this difference vanishes
fast with $\lambda$ and one is not able to almost surely detect a difference
between the two macroscopic interfaces, or the difference between these two
behaviors can be detected by averaging them out over the $\lambda^{7 / 4}$
parts of the path. In the end, one has to compare certain critical exponents
to decide which scenario is correct and it turns out that it is  the second one.
A ``flavor'' of supercritical percolation is therefore still present in the
scaling limit.

We would like to emphasize that we do not really use any complex analytic SLE technology
to derive our results. Our considerations are based on percolation techniques, 
Kesten's scaling ideas and the knowledge of the exact value of the critical
exponents (the derivation of which however used SLE).

\section{Preliminaries}
\label {prel}

This paper will build on earlier results on near-critical percolation, and in particular on Kesten's paper \cite {Ke4}. In this section, we recall some results that we shall use. All these results are stated and derived in \cite {N2}. See also the last section of \cite {Wpc} for a rough survey. Throughout the paper, we will use the notation $\asymp$ to say that the ratio between two quantities remains bounded away from $0$ and $\infty$.

Recall our definition of the correlation length
\[ L (p) = L (p, \epsilon) = \inf \{N : R (p, N) > 1 / 2 + \epsilon\}. \]
Note that $L$ is non-increasing with respect to $p$ and non-decreasing with respect to $\epsilon$.

We define the event $A^4 (n)$ that there exist four paths of alternating colors
(when ordered clockwise around the origin), each of them joining a neighbor of the origin
to the circle of radius $n$. Similarly, we define the event $A^2 (n)$ (resp. $A^1 (n)$) that there
exist two paths of different colors (resp. one black path) from a neighbor of the origin to the circle
of radius $n$. We will use the following properties, that hold for any fixed $\epsilon \in (0,1/2)$:
\begin{enumerate}

  \item The Russo-Seymour-Welsh estimates remain valid below $L(p,\epsilon)$: for all $k \geq 1$, there exists $\delta_k(\epsilon) >0$ such that for all $p$, for all $N \le L (p,\epsilon)$, the probability (for $P_p$) that there is a black (resp. white) horizontal crossing of a $k N \times N$ parallelogram is at least $\delta_k$.

  \item Let $A^2 (n_1, n_2)$ denote the event that there exist two arms of
  different colors joining the circles of radii $n_1$ and $n_2$. Then (this follows from the ``arm-separation lemmas''),
  \[ P_p (A^2 (n_1 / 2)) \times P_p (A^2 (2 n_1, n_2)) \asymp P_p (A^2 (n_2)) \]
  uniformly for $p \geq 1/2$ and $2 n_1 \le n_2 \le L(p,\epsilon)$ (this is known as the \emph{quasi-multiplicativity} property). The same is true for four arms and one arm, \emph{i.e.} with $A^4$ or $A^1$ instead of $A^2$. \label{multip}

  \item It can be used to prove that
  $$\sum_{j=1}^n j P_p ( A^2 (j)) \asymp n^2 P_p ( A^2 (n)),$$
  uniformly for $p \ge 1/2$ and $n \le L(p,\epsilon)$. \label{summation}
  
  \item We have
  \[ P_p (A^2 (n)) \asymp P_{1 / 2} (A^2 (n)) \text{ and } P_{p} (A^4 (n))
     \asymp P_{1 / 2} (A^4 (n)) \]
  uniformly for $p \ge 1/2$ and $n \le L (p, \epsilon)$. \label{arms_near}

  \item Finally, $$(p - 1 / 2) L (p, \epsilon)^2 P_{1 / 2} (A^4 (L (p,
  \epsilon))) \asymp 1$$ as $p \to 1/2^+$. \label{equiv}

\end{enumerate}

The last property, which holds for each $\eps \in (0,1/2)$, implies in particular that for any fixed $\eps,\eps' \in (0,1/2)$,
$$p^* (N, \eps) -1/2 \asymp p^* (N, \eps') - 1/2$$
as $N \to \infty$, and that
 $$L (p, \epsilon) \asymp L (p, \epsilon')$$ as $p \to 1/2^+$.
 Note that the combination of items \ref{multip}. and \ref{arms_near}. shows that 
$$P_p (A^2 (n_1, n_2)) \asymp P_{1/2} (A^2 (n_1,n_2))$$
uniformly for $p \ge 1/2$, $2n_1 \le n_2 \le L(p, \eps)$ (and similarly for $A^4$).  

Recall also {\cite{SmW}} (see also {\cite{Wpc}}) that for any $\eta \in (0,1)$,
\[ P_{1 / 2} (A^2 (\eta n,n)) \to f_2 (\eta)  \text{ and } P_{1 / 2} (A^4 (\eta n,n)) \to f_4 (\eta)\]
as $n \to \infty$, where $f_2 (\eta ) \sim \eta^{1/2+ o(1)}$ and $f_4 (\eta) \sim \eta^{5/4+ o(1)}$ as $\eta \to 0^+$,  which in turn implies that
\[ P_{1 / 2} (A^2 (n)) = n^{- 1 / 4 + o (1)} \text{ and } P_{1 / 2} (A^4 (n)) =
   n^{- 5 / 4 + o (1)} \]
as $n \to \infty$.
Analogous statements hold for any number of arms (see {\cite{N2}}).

Note that the correlation length introduced here differs slightly from the usual one -- let us denote it here by $L^{\ast}(p,\epsilon)$ -- that is defined in terms of crossing probabilities of rectangles or rhombi (the one introduced in \cite{CCF}, and used for instance in \cite{Ke4}). Using Russo-Seymour-Welsh considerations (and also the fact that  for any $\eps,\eps' \in (0,1/2)$, $L^{\ast}(p,\epsilon) \asymp L^{\ast}(p,\epsilon')$ -- see \cite{Ke4,N2}), these two definitions can easily be shown to be equivalent: for any $\eps \in (0,1/2)$, $L(p,\epsilon) \asymp L^{\ast}(p,\epsilon)$. We are thus allowed to use for $L$ results established for $L^{\ast}$.

Alternatively, we could have chosen to work with rhombi instead of triangles and then, we would have used directly the usual definition of correlation length. These two definitions are also known to be equivalent to other ``natural'' correlation lengths, describing for example the mean radius of a finite cluster or the rate of decay of connectivity properties (see e.g. the discussion in section 2.2 of \cite{CNo}).

\section{Tightness}

For simplicity of presentation, we will stick to the setup that we described in the
introduction, even if the shape of the considered domain (in our case, it is the triangle) could be chosen
arbitrarily. 

Let us describe this setup more precisely. We follow here the
definitions of {\cite{AB}}, and refer the reader to this paper for more
details. We consider $T$ the (filled) equilateral triangle of unit side length, with
corners $(1 / 2, 0)$, $(- 1 / 2, 0)$ and $(0, \sqrt{3} / 2)$. The rescaled
interfaces will be elements of $\mathcal{S}_T$, the space of curves in $T$:
these are equivalence classes of continuous functions from $[0, 1]$ to $T$,
where two functions $f_1$ and $f_2$ represent the same curve if and only if
there exists a continuous increasing bijection $\phi : [0, 1] \to
[0, 1]$ such that $f_1 = f_2 \circ \phi$. We endow this space
with the quotient metric:
\[ d (f_1, f_2) :=  \inf_\phi \bigg( \max_{[0,1]} |f_1 - f_2 | \bigg) \]
where the infimum is over the set of increasing bijections $\phi$ from $[0,1]$ onto itself.
 
We call $P_{p, N}$ the law of the rescaled interface $\Gamma_N : =
\gamma^N / N$ of percolation with parameter $p$ in our triangle (with mesh size $1 / N$): this is a probability measure on $\mathcal{S}_T$. Endowed
with the previous metric, $\mathcal{S}_T$ is a complete separable space, so
that tightness and relative compactness are equivalent (by Prohorov's
theorem).

\begin {prop}
The family $(P_{p,N}, p \in [0,1], N \ge 1)$ is relatively compact in ${\mathcal S}_T$.
\end {prop}

\noindent
{\bf Proof.}
Let us consider a sequence $(P_{p_k,N_k})$ in
this family. Our goal is to find a converging subsequence:
\begin{itemize}
  \item By symmetry, we can assume that $p_k \ge 1 / 2$ for all $k \ge 1$.
  
  \item If $N_k$ remains bounded along a subsequence, proving convergence of a subsequence of
  $(P_{p_k, N_k})$ is trivial.
  
  \item We can therefore restrict ourselves to the case where $N_{k} \to
  \infty$. Suppose that for all $\epsilon \in (0,1/2)$, $N_k / L (p_k, \epsilon) \to \infty$ along a
  subsequence; then it is easy to check (using Russo-Seymour-Welsh arguments)
  that along that subsequence, $\Gamma^k= \Gamma_{N_k}$ converges in law to the
  concatenation of the two parametrized segments $[0, 1 / 2] \cup [1 / 2,
  \sqrt{3} i / 2]$.
  
  \item We now suppose that for some $\epsilon \in (0, 1/2)$, $N_{k} / L (p_{k},
  \epsilon)$ remains bounded. In particular, we get that for some fixed $\epsilon' > 0$, $N_{k} \le L (p_{k},\epsilon')$ for all large enough $k$ (indeed, by Russo-Seymour-Welsh, one can increase $L$ by at least any constant factor by choosing  a larger $\epsilon'$).
  We can then use the machinery developed by 
  Aizenman and Burchard (Theorem 1.2 in \cite {AB}): 
  Russo-Seymour-Welsh estimates hold uniformly for all $p$ and $N \le L (p,
  \epsilon')$, so that for any annulus 
  ${\cal A} (x; r, R) =\{z \in T : r <
  |z - x| < R\}$,
  \[ P_{p_{k}} ( {\cal A} (x; r, R)  \hbox { is traversed by }\Gamma_{N_{k}})
  \le  P_{p_k} ( A^1 (rN_k, RN_k  ) ) \le
      C (r / R)^{\alpha} \]
  for two universal constants $\alpha, C$ (traversed means here that it visits the inner and the outer boundary of the annulus). The BK inequality (see e.g.
 {\cite{Grbook}}) then leads to the fact that
  \[ P_{p_{k}} ( {\cal A}  (x; r, R) \hbox { is traversed } K 
  \hbox { times by
     } \Gamma_{N_{k}} ) \leq C_K (r / R)^{\alpha K} \]
  which is exactly the hypothesis (H1) of {\cite{AB}} (uniform power bounds on
  the probability of multiple crossings in annuli).
  
  For the sake of completeness, let us just sketch why tightness can be
  derived from this. First, the hypothesis (H1) implies regularity properties
  (Theorem 1.1 in \cite {AB}) for the random curves $\Gamma_{N_{k}}$: for any $\beta>0$, there is some H\"older continuity bound such that for each $k$, there is a probability at least $1-\beta$ that $\Gamma_{N_{k}}$ can be parametrized and satisfy this bound. In particular, we can exhibit some equicontinuous set such that for each $k$, $\Gamma_{N_{k}}$ belongs to this set with probability at least $1-\beta$. Tightness follows (Theorem 1.2 in \cite {AB}) using (a slight adaptation of) Arzela-Ascoli's characterization of compactness for continuous functions on a compact set.
\end{itemize}
Hence, in all cases, one can find a converging subsequence of $(P_{p_{k}, N_{k}})$.
\qed

\section{Length and dimension of near-critical interfaces}

 In this section, we will study the non-degenerate case where $p(N)- 1/2$ is not too large: we assume that $p(N) \in [1/2, p^* (N, \eps)]$ for some fixed $\eps \in (0,1/2)$ (note that this includes the critical case). We exhibit properties that near-critical interfaces share with critical ones.
We will first focus on an estimate for discrete interfaces and then study the scaling limits.

\subsection{Length of discrete interfaces}

The length (number of edges) $\ell (\gamma^N)$ of the discrete interface $\gamma^N$ in the finite-size
scaling regime is roughly the same as that of the critical interface:

\begin {prop}
\label{proplength}
Assume that $p(N) \in [1/2, p^* (N, \eps)]$. Then for any fixed $\beta >0$,
$$\lim_{N \to \infty } P_{p(N)} ( \ell (\gamma^N) \notin [N^{7/4 - \beta}, N^{7/4 + \beta }]) = 0.$$
\end {prop}

\medbreak
\noindent
{\bf Proof.}
Note that the hypothesis on $p(N)$ implies in particular that $N \le L(p(N),\eps)$, so that we can use the results on near-critical percolation that we have recalled in the preliminaries.

\begin{figure}
\begin{center}
\includegraphics[width=8cm]{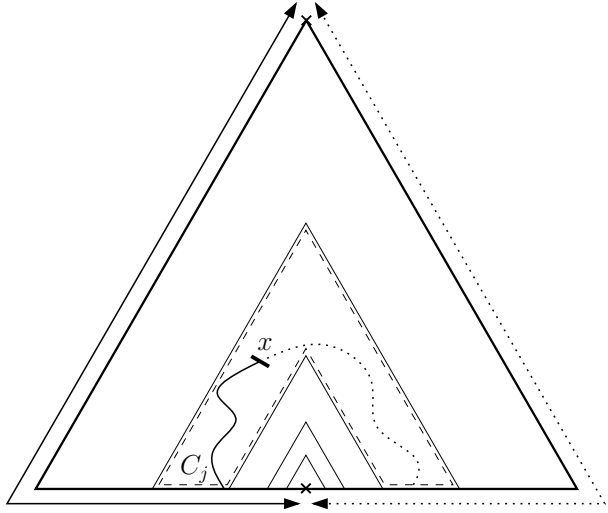}
\caption{For an edge $x$, being on $\gamma^N$ corresponds to the existence of two arms. For the lower bound, we consider semi-annuli $C_j = T_{3^j} \setminus T_{3^{j - 1}}$ and we restrict to the edges at distance at least $3^{j - 3}$ from the boundary of $C_j$ (dashed).}
\label{2arms}
\end{center}
\end{figure}

Let us first derive the upper bound for $\ell( \gamma^N) $: we estimate its expectation. We note that for any edge $x \in T_N$ on the dual hexagonal lattice, $x \in \gamma^N$ if and only if  there exist two arms as depicted on Figure \ref{2arms}, one black joining a neighbor of $x$ to the negative real half-axis, and one white joining a neighbor of $x$ to the positive real half-axis.

In particular,
\[ P_p (x \in \gamma^N) \le P_p (A^2 (d (x, \partial T_N))) . \]
It allows us to get an upper bound on $E_p(\ell( \gamma^N)) $: for each $\beta> 0$, and for some constants $c, c'\ldots$ (using the estimates on near-critical percolation and the fact that there are never more than $3N$ points at distance $j$ from $\partial T_N$),
\begin{eqnarray*}
  E_p (\ell (\gamma^N)) & \le & \sum_{x \in T_N} P_p (A^2 (d (x, \partial
  T_N)))\\
  & \le & c \sum_{x \in T_N} P_{1/2} (A^2 (d (x, \partial T_N)))\\
  & \le & c' \sum_{j = 1}^N N  P_{1/2} (A^2 (j))\\
  & \le & c'' N \sum_{j = 1}^N j^{ - 1 / 4 + \beta/2 }\\
  & \le & c''' N^{7 / 4 + \beta/2} .
\end{eqnarray*}
The upper bound for $\ell( \gamma^N) $ follows immediately by Markov's inequality.

\medbreak

For the lower bound, we use the standard second moment method. For that purpose, take $l$ such that $3^l \le N < 3^{l+1}$ and decompose the triangle $T_N$ into concentric ``semi-annuli'' $C_j := T_{3^j} \setminus T_{3^{j - 1}}$, $j \leq l$. Note that the curve $\gamma^N$ has to cross each of these semi-annuli.

For any edge $x$ in $C_j$ at distance at least $3^{j - 3}$ from the
boundary of $C_j$, the results that we have recalled together with the
``separation lemmas'' imply that the probability $Q (x)$ that there exist two
paths as before, from $\partial x$ to the positive and negative real half-axes, that stay in $C_j$, satisfies
\[ Q (x) \asymp P_{1 / 2} (A^2 (3^j)) \]
uniformly for such $x$ and $3^j \le L (p)$. Such an edge  $x$ is
necessarily on $\gamma^N$ by the previous remark. Let $\mathcal{N}_j$ denote the number of these
edges: we have 
$$E_p ( \mathcal{N}_j) \asymp (3^j)^2 P_{1 / 2} (A^2 (3^j)) \ge (3^j)^{7/4 - \beta}$$ 
for large $j$ (and any fixed $\beta>0$).

Now, using the quasi-multiplicativity property and then item \ref{summation}. we get
\begin{eqnarray*}
  E_p ( \mathcal{N}_j^2) & \le & \sum_{x, x'} P_p (A^2 (d (x,
  x') / 2))^2 P_p (A^2 (2 d (x, x'), 3^{j - 3}))\\
  & \le & c E_p ( \mathcal{N}_j) \times \sum_{k = 1}^{3^j} k P_p (A^2
  (k))\\
  & \le & c' E_p ( \mathcal{N}_j)^2 .
\end{eqnarray*}
It follows that for some (universal) constant $c_0>0$,
\[ P_p ( \mathcal{N}_j \ge c_0 E ( \mathcal{N}_j)) > c_0 . \]
Consider some $j_0$: since the events $\{ \mathcal{N}_{l-j} \ge c_0 E_p ( \mathcal{N}_{l-j} ) \}$ for $ j \in \{0, \ldots, j_0 \}$ are independent, we get that
$$ P_{p(N)} ( \mathcal{N}_{l-j} \le c_0 E_p ({\cal N}_{l-j}) \text{ for all $j \le j_0$}) \le (1- c_0)^{j_0+1}.$$
The lower bound follows readily because 
$$ \ell (\gamma^N) \ge \max ( {\cal N}_l, {\cal N}_{l-1}, \ldots )$$
and our lower bound for $E ( {\cal N}_j )$.
\qed

\subsection{Dimension of scaling limits}

We now show that scaling limits of near-critical interfaces have the same Hausdorff dimension as in the critical regime:

\begin {prop}
Assume that the law of the curve $\gamma$ is the weak limit of $(P_{p_k,N_k})$, with $N_k \to \infty$ and $p_k \in [1/2, p^* (N_k, \eps)]$. Then the Hausdorff dimension of $\gamma$ is almost surely $7/4$.
\end {prop}

Similar statements for $n$-point correlation functions, multiple crossings, multifractral spectrum etc. could be derived in the same way. We leave this to the interested reader (see the informal discussion at the end of the paper) and we focus here only on this fractal dimension.

We might remark that the a priori estimates for the existence of two or three 
arms near a half-plane, or for the existence of five arms in the interior of the domain still hold (uniformly) for $n \le L(p, \eps)$,
 as these are consequences of the Russo-Seymour-Welsh estimates only. 
Hence, exactly as in the convergence to SLE(6) case (e.g. \cite {Sm2,Wpc}), the discrete hitting probabilities converge to the continuous hitting probabilities.
  
\medbreak
\noindent 
{\bf Proof.}
The argument goes along similar lines as for the discrete length.
For the upper bound, it suffices to prove that for all  $\beta>0$,  the expected number of balls of radius $\eta$ needed to cover $\gamma$ is bounded by $\eta^{-7/4 - \beta}$. But, we know that this expected number is bounded uniformly by this quantity for the discrete paths $\Gamma_{N_k}$, and so this also holds in the scaling limit.
   
For the lower bound, we decompose the triangle into concentric semi-annuli, and we exhibit a family $({\cal C}_j, j \ge 1)$ of independent 
events of probability at least $c_0$, such that on each ${\cal C}_j$, the curve $\gamma$ has dimension $7/4$.
Again, we follow exactly the same idea (with second moment estimates) as in the discrete case, using the fact (noticed in Section \ref {prel}) that the discrete estimates for $P_p (A^2 (n_1, n_2))$ hold uniformly. 
\qed

\section {The alternative}

We now prove that for the scaling limits of percolation interfaces, there are only three main options. 
\begin {prop}
\label {dichotomy}
Suppose that $P$ is the weak limit of a sequence 
 $(P_{p_k, N_k})$ when $k \to \infty$, with $N_k \to \infty$ and $p_k \geq 1/2$. 
Then:
\begin{itemize}
  \item either $P$ is the Dirac mass on the concatenation of the two
  boundary segments $[0,1/2]$ and $[1/2, i \sqrt {3}/2]$,
  \item or $P$ is the law of SLE(6),
  \item or for some $0 < \epsilon' < \epsilon < 1/2$, 
  $p(N_k, \epsilon')  \le p_k  \le p( N_k, \epsilon)$ for all large $k$.
\end{itemize}
\end {prop}

\noindent
{\bf Proof.}  
Suppose now that $P$ is the limiting law of a sequence $\Gamma^k:= \Gamma_{N_k}$
defined under $P_{p_k, N_k}$, where $N_k \to \infty$ and $p_k \ge 1/2$. 
\begin {itemize}
\item
Let us first suppose that for each $\epsilon \in (0,1/2)$, $N_k \ge L( p_k, \epsilon)$ for all sufficiently large $k$. This implies readily (by standard Russo-Seymour-Welsh arguments) that for each fixed $\epsilon'$, $N_k / L (p_k, \epsilon') \to \infty$. It follows readily that the probability that there exists a closed cluster of radius greater than any fixed positive $u$ for $P_{p_k, N_k}$ in the rescaled triangle goes to $0$ as $k \to \infty$ and the $P$ must be supported on the path that goes along the right boundary of the triangle.
\item
Let us now suppose that for each $\epsilon >0$, $N_k \le L(p_k, \epsilon)$ for all sufficiently large $k$.  
 We assume that the law $P$ is not the law of an SLE(6) and try to reach a contradiction. 
 In the discrete setting, it is always possible to couple $\Gamma^k$ (corresponding to $P_{p_{k}}$)
 with a realization of the interface $\tilde \Gamma^k$ for $P_{1/2}$ in such a way that the sites that are black for the realization of $P_{1/2}$ are also black for the realization of $P_{p_k}$. In this way, the path $\Gamma^k$ is more ``to the right'' than $\tilde \Gamma^k$. Letting $k \to \infty$, we see that we can couple a realization $\Gamma$ of $P$ with an SLE(6) $\tilde \Gamma$ in the triangle $T$ with the same property.
 
It follows from our results on Hausdorff dimensions (in fact the uniform Russo-Seymour-Welsh bounds would suffice) that $\Gamma$ is almost surely  a continuous curve with zero Lebesgue measure. Clearly, one can divide the set of connected components of $T \setminus \Gamma$ into two parts. Intuitively speaking, those ``to the right of $\Gamma$'' and those ``to the left of $\Gamma$''. Let us call ${\cal R} (\Gamma)$ and ${\cal L} (\Gamma)$ the respective unions of these components. 
Similarly, we can define ${\cal R} (\tilde \Gamma)$ and ${\cal L} (\tilde \Gamma)$. Then, we have that ${\cal L} (\tilde \Gamma ) \subset 
{\cal L} (\Gamma)$.

The fact that the laws of $\Gamma$ and $\tilde \Gamma$ are not identical implies that with positive probability the two open sets 
${\cal L} (\tilde \Gamma)$ and ${\cal L} (\Gamma)$ are not equal. 
It follows (using also the fact that $\Gamma$ is a curve that does not backtrack on its own past) 
that for some $z \in T$ with rational coordinates and for some $r >0 $ with rational coordinates, if we denote by $D(z,r)$ the disc centered at $z$ with radius $r$, 
$$ P ( D(z,r) \subset {\cal L} (\Gamma) \hbox { and } D ( z, r) \subset {\cal R} (\tilde \Gamma) ) > 0 .$$
We also choose $r$ small enough so that the distance between $z$ and the boundary of the triangle is at least $4r$.
From now on, $z$ and $r$ will be fixed.

The convergence of the discrete interfaces to the continuous ones imply that for some positive constant $c$, 
$$ P ( D(z,r/2) \subset {\cal L} (\Gamma^k) \hbox { and } D ( z, r/2) \subset {\cal R} (\tilde \Gamma^k) ) > c $$
for any large $k$, with obvious notation. Note that these discrete events 
do not depend on the state of the sites inside $D(z,r/2)$ as they imply that neither $\Gamma^k$ nor $\tilde \Gamma^k$ did intersect the disc. 
So, it follows from the uniform RSW estimates that if we define the event ${\cal H}^k$ that there exist a white circuit in $D(z,r/2) \setminus D(z, r/4)$ for $P_{p_k}$ and a black circuit for $P_{1/2}$ in the same annulus (where both circuits  disconnect $z$ from $\partial D(z,r/2)$), we get that for some absolute positive constant $c'$ and all sufficiently large $k$,
$$ P ( {\cal H}^k \hbox { and } D(z,r/2) \subset {\cal L} (\Gamma^k) \hbox { and } D ( z, r/2) \subset {\cal R} (\tilde \Gamma^k) ) > c' .$$

It follows readily that if one increases continuously $p$ from $1/2$ to $p_k$, there will exist -- with probability bounded from below -- a value $p$ at which the interface jumps from the left to the right of $D(z,r/2)$ without hitting it. At this $p$, there exist necessarily four arms of alternating colors originating at the site $x$ that has been flipped. The different possibilities are depicted in Figure \ref{flip}. 

\begin{figure}[t]
\begin{center}
\includegraphics[height=4.3in]{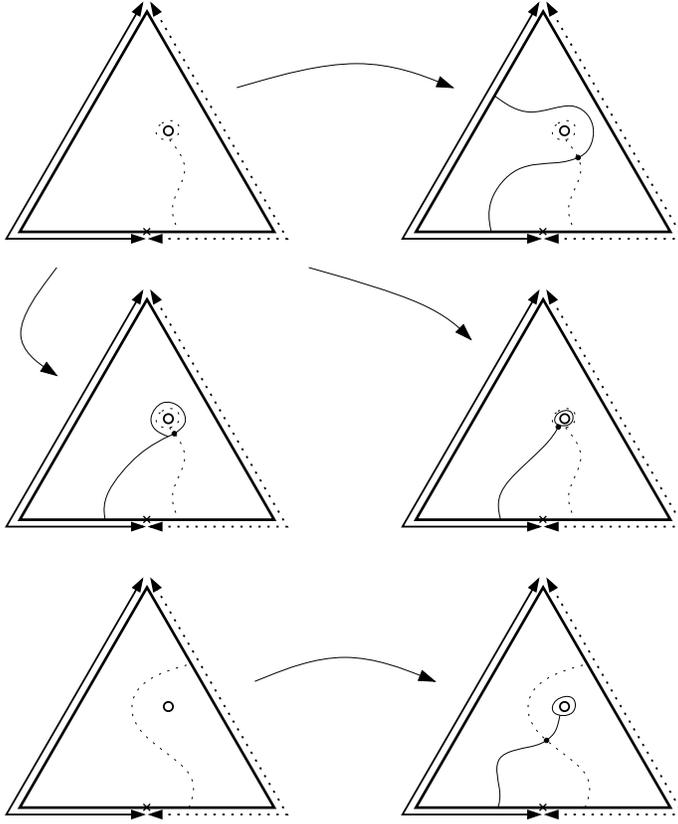}
\caption{Different possibilities to jump from the left to the right of the small disc $D(z,r/2)$.}
\label{flip}
\end{center}
\end{figure}

For each given site $x$ on the grid $T_{N_k}/N_k$, let us evaluate the probability that this occurs.
 When $x$ is not close to the boundary of the triangle, i.e.  if it is at distance at least $r/2$, then we see that four arms of alternating colors must originate from $x$ all the way to distance $r/2$ from $x$. The probability is then bounded by 
$P_{p} ( A^4 (2N_k / r)) \asymp P_{1/2} (A^4 (N_k))$ (recall that $r$ is fixed and that $p \le p_k$).  

When $x$ is at distance $j/N_k < r/2$ from the boundary, then one must have four arms of alternating colors up to distance $j/N_k$, and then two arms that cross a semi-annulus (in a half-plane) from radius $j/N_k$ to $r/2$. Combining this with quasi-multiplicativity and the known bounds (2 arms in a half-plane, 5 arms in the plane that give an upper-bound for the four-arm probability) leads to a bound of the type
$$ c P_{1/2} (A^4 (N_k)) \times (j/rN_k) \times (j/N_k)^{-2+\beta}.$$

For any $j$, there are at most a constant times $N_k$ points on our grid that are at distance exactly $j/N_k$ from the boundary of the triangle. 
Hence, when we sum over all possible sites $x$ in the triangle, we get that the probability that between $p$ and $p+dp$ a site $x$ is flipped in such a way that the curve jumps from the left to the right of the small disc is bounded by $dp$ times
$$
c' P_{1/2} ( A^4 ( N_k))\bigg( N_k^2  + N_k \times N_k^{1 - \beta} \sum_{j=1}^{N_k} j^{-1 + \beta} \bigg)
\le c'' N_k^2 P_{1/2} ( A^4 ( N_k))$$
(for some constants $c', c''$ that do not depend on $k$).

It finally follows that for some positive constant $c$,
$$
(p_k - 1/2) \times N_k^2 \times P_{1/2} ( A^4 ( N_k)) \ge c.
$$ 
But if we apply Russo's formula in a similar way to evaluate $R(p_k, N_k) - R(1/2, N_k)$ together with the fact that 
the four-arm probabilities are uniformly controlled for $n \le L(p)$, we get that for some  positive $c'''$,
$$
R(p_k, N_k) - R( 1/2, N_k) \ge c''' ( p_k -1/2) \times N_k^2 \times P_{1/2}(A^4 (N_k)).
$$
Hence, there exists a positive constant $c''''$ such that 
$R(p_k, N_k) \ge 1/2 + c''''$ for all sufficiently large $k$. In other words, 
$N_k \ge L(p_k, c'''')$ for all sufficiently large $k$, and this contradicts our assumption.

\item
Finally to conclude the proof of  the proposition, let us suppose that we are not in the last scenario described in its statement i.e. that it is not true that for some $0 < \eps' < \eps < 1/2$, 
$ p(N_k, \eps') \le p_k \le p(N_k, \eps)$ for all large $k$. Then it means that a subsequence of $(p_k, N_k)$ falls into one of the two cases that we have just studied. Hence, the limiting probability measure  is either the law of SLE(6) or the Dirac mass on the union of the two boundary segments.
\end {itemize}
\qed

\section{Asymmetry of near-critical interfaces}

In this section, we are going to assume that we are in the intermediate regime where the discrete percolation 
parameter $p$ used at scale N satisfies $$p(N) \in [p^* (N,\eps), p^* (N,\eps')]$$
for some fixed $0 < \eps < \eps' < 1/2$. The previous alternative tells us that these are the ``non-trivial'' cases in the large-scale limit. 
We are going to exhibit asymmetry properties that enable to 
detect the difference with a critical interface. We will first derive a result in the discrete setting, and then a result in the continuous setting for the scaling limit.

\subsection{Discrete asymmetry}

On the discrete level, it is possible to explore the interface dynamically:
each time one discovers a new hexagon, one tosses a coin to decide if it is black
or white, so that with probability $p$ (resp. $1 - p$), one makes a 60 degree
turn to the right (resp. to the left). Of course, one often bumps into an
already discovered hexagon, and in this case, the turn is already determined by a previous tossing. Hence, the percolation interface is the image
under some map of a sequence of coin tosses. 
Let us call $\ell^+ (\gamma^N)$ and $\ell^- (\gamma^N)$ the number of black
(resp. white) cells neighboring the path $\gamma^N$.

\begin {prop}
\label {das}
Suppose that $p(N) \in [p^* (N,\eps), p^* (N,\eps')]$. Then for any fixed $\beta>0$,
$$\lim_{N \to \infty} P_{p(N)} ( \ell^+ (\gamma^N) - \ell^- (\gamma^N) \ge N^{1-\beta} ) = 1 .$$
On the other hand, for critical percolation, 
$$
\lim_{N \to \infty} P_{1/2} (  \ell^+ (\gamma^N) - \ell^- (\gamma^N)  \le N^{7/8 + \beta }) = 1.
$$\end {prop}

\noindent
{\bf Proof.} Define $v(N) = p(N) - 1 / 2$. Clearly the difference between the number of black 
neighbors discovered by $\gamma^N$ and the number of white discovered neighbors evolves like a biased
simple random walk with drift $2 v$ (at each step the probability to add
one is $1 / 2 + v$) stopped after a certain random number of steps $\ell$
that is roughly of the order $N^{7 / 4}$ (see Proposition \ref{proplength}).
Hence, we see that this difference will be of the order
\[ 2 v \times \ell + O ( \sqrt{\ell}) \]
(with high probability) i.e. more precisely that for all $\beta >0$,
$$\lim_{N \to \infty} P_{p(N)} ( \ell^+ (\gamma^N) - \ell^- (\gamma^N) \ge 2 v \ell - \ell^{1/2+ \beta/2}  ) = 1 .$$
In the off-critical case where $v(N) \ge N^{- 3 / 4 - \beta/2}$,
we see  that the drift term dominates because of Proposition \ref {proplength} so that the first statement follows readily.
Similarly, when $v=0$, the result follows directly from Proposition \ref {proplength}.
\qed

\medbreak
Hence, the discrete interface in the off-critical regime where $p(N) \ge
p^{\ast} (N,\epsilon)$ will have (on a microscopical scale) anomalously more
``black neighbors'' than ``white neighbors''. Moreover, the proof shows that one can
detect a good approximation of the actual value of $p(N)$ by looking at this
difference because  $(\ell^+ (\gamma^N) - \ell^- (\gamma^N)) / \ell$ 
will be close to $v(N)$.

But this result is not sufficient to deduce asymmetry  of  scaling limits.
Indeed, if we had for instance taken
\[ p(N) = 1 / 2 + N^{- b} \]
for some $b \in (3 / 4, 7 / 8)$, then an analogous discrete asymmetry result would
hold i.e. for some small $\beta >0$
$$
\lim_{N \to \infty} P_{p(N)} ( \ell^+ (\gamma^N) - \ell^- (\gamma^N) \ge N^{7/8 + \beta}  ) = 1,
$$
but we have proved in the previous section that in this case, the interface still converges to SLE(6) as $N \to \infty$ (and we know that SLE(6) is symmetric).

\subsection{Continuous asymmetry}

The goal of this section is to prove the following result:

\begin {prop}
Suppose that $\gamma$ is the limit (in law) of a sequence of $\Gamma_{N_k}$'s where $p_k \in [p^* (N_k,\eps), p^* (N_k,\eps')]$ (and $N_k \to \infty$).
Then the law of $\gamma$ is singular with respect to that of SLE(6).
\end {prop}

\begin{figure}
\begin{center}
\includegraphics[width=12cm]{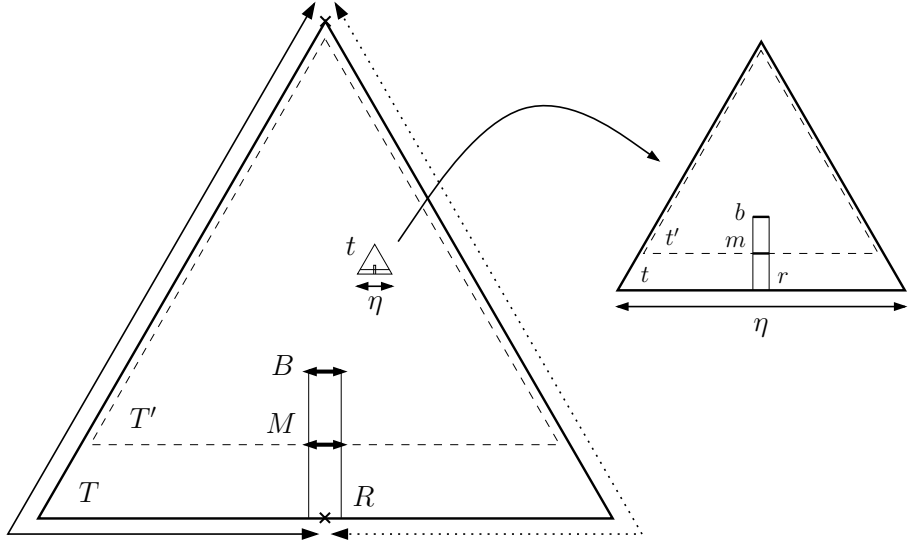}
\caption{In the triangle $T$, we consider the subsets $R$, $M$, $B$ and $T'$. In any small equilateral triangle $t$ of size $\eta$, we define similarly the corresponding sets $r$, $m$, $b$ and $t'$.}
\label{subsets}
\end{center}
\end{figure}

\noindent 
{\bf Proof.}
Let us fix a very small positive $\beta$.
Define the rectangle $R=[-a, a]\times [0, 1/4]$ for a small fixed $a$. Define the segments $M=[-a, a] \times \{1/8\}$ and $B= [-a, a]\times \{1/4\}$, and the smaller triangle $T'= \{ (x,y) \in T \ : \ y \ge 1/8 \}$ (see Figure \ref{subsets}). For each small equilateral triangle $t$ of size $\eta$ (i.e. that one can obtain by translating $\eta T$),  we define analogously the sets $r$, $m$, $b$ and $t'$. We say that the triangle $t$ is \emph{good} for $\gamma$, if at the first hitting time $\sigma(t)$ of $t \setminus r$, $\gamma (\sigma(t)) \in b$, and if 
$\gamma$ does exit $r$ through $b$ after its first hitting time of $m$.

\begin{figure}
\begin{center}
\includegraphics[width=9cm]{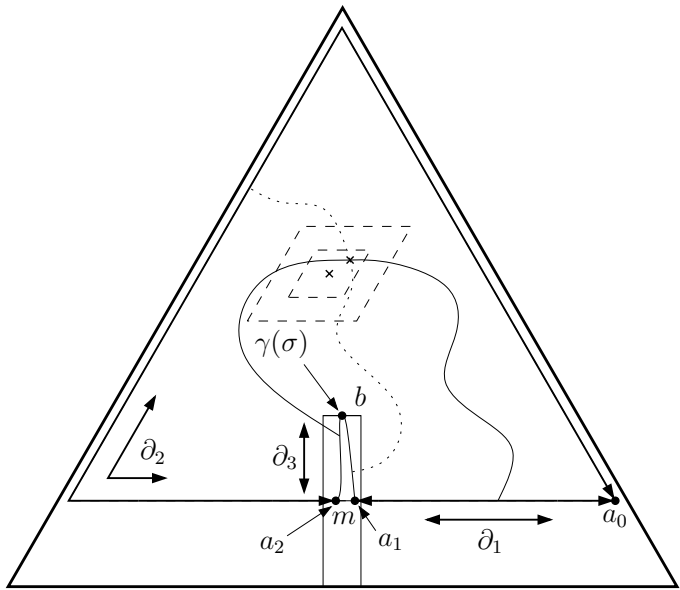}
\caption{Using ``pivotal'' sites like the one depicted here provides a lower bound on the conditional probability that the small triangle $t$ is very good.}
\label{pivotal}
\end{center}
\end{figure}

Let us define the connected component $d$ of $t' \setminus \gamma [0, \sigma]$ that has the top boundaries of $t'$ on its boundary. We consider the following three marked points on the boundary of $d$: suppose that we go along the boundary of $d$ clockwise starting from the bottom right corner of $t'$ (that we call $a_0$). Then $a_1$ is the right-most  point on $m \cap \gamma (0,\sigma)$ and $a_2$ is the left-most point on $m \cap \gamma (0,\sigma)$ (see Figure \ref{pivotal}). We call $\partial_1$ (resp. $\partial_2$) the part of the boundary of $d$ between $a_0$ and $a_1$ (resp. between $a_2$ and $a_0$).

When $t$ is a good triangle, we define $F(\gamma,t)$ to be the probability given by Cardy's formula of the existence in $d$ of a crossing from $\partial_1$ to the part $\partial_3$ of the boundary between $\gamma (\sigma) $ and $a_2$. 
If after $\sigma$, $\gamma$ does hit $\partial_1$ before $\partial_2$, we then say that $t$ is \emph{very good} for $\gamma$.
If $\gamma$ was an SLE(6), then $F(\gamma,t)$ would just be the conditional probability given $\gamma ([0, \sigma])$ that $t$ is very good for $\gamma$.

But here, the law of $\gamma$ is the limit of the law of near-critical percolation interfaces $\Gamma^k= \Gamma_{N_k}$ defined under $P_{p_k, N_k}$, where $p_k \ge p^* ( N_k, \eps)$ and $N_k \to \infty$. 
By Skorokhod's representation theorem, we can couple on the same probability space samples of all  $\Gamma^k$'s  and $\gamma$ 
in such a way that $\Gamma^k \to \gamma$ almost surely just as in the proof of convergence of critical percolation interfaces to SLE(6) (as in \cite {Wpc} for instance). For $\Gamma^k$, we can also check if a triangle $t$ is good and very good. In fact, for each fixed triangle, the probability that its status for $\gamma$ is not identical to its status for $\Gamma^k$ goes to zero as $k \to \infty$. Indeed, if we suppose for instance that $t$ is good for $\gamma$ but not for $\Gamma^k$, then it means that a three-arm event occurs for $\gamma$ and a priori bounds do exclude this. 

If we define the stopping time $\sigma^k$ associated to the triangle $t$ and the path $\Gamma^k$, this means that when $t$ is good, then $\Gamma^k [0, \sigma^k]$ converges almost surely to $\gamma [0, \sigma]$. So, the domains $d^k$ (with obvious notation) converge in Carath\'eodory topology towards $d$ almost surely.

We now couple each realization of $\Gamma^k$ with a realization $w^k$ of critical percolation on the same lattice, in such a way that the realization used to define $\Gamma^k$ dominates $w^k$
 (\emph{i.e.} its set of black sites is larger). Then (see e.g. \cite {Wpc}), it follows that the probability $F^k (\Gamma^k, t)$ that there exists a crossing from $\partial_1^k$ to $\partial_3^k$ in $d^k$ for $w^k$ 
(with natural notation) converges almost surely to $F(\gamma,t)$.

Suppose now that the triangle $t$ is good for $\Gamma^k$ and consider the configurations as depicted on Figure \ref{pivotal}: we restrict ourselves to the case where there exist sites in a smaller rhombus of size comparable to $\eta N_k$ from which four arms of alternating colors originate. For each such site $x$, and conditionally on $\Gamma^k [0, \sigma^k]$, the probability that this four-arm configuration exists for percolation of parameter $p \in [1/2, p_k]$ in the domain not yet explored by $\Gamma^k$ is bounded from below by a constant $c$ times 
$P_{1/2} (A^4 ( \eta N_k))$, uniformly over $p$, $k$ and the choice of $p$ (this follows from the uniform estimates for four arms, the separation lemmas and Russo-Seymour-Welsh).
The probability that this site is flipped from white to black when one increases $p$ from $1/2$ to $p_k$, and that at this value $p$ the configuration is as described above, is 
therefore larger than 
$$ (p_k -1/2) \times c \times P_{1/2} (A^4 (\eta N_k)).$$
Finally, note that this event can happen for only one site (when $p$ increases, since the crossing event is increasing, there can exist only one value of $p$ and therefore one site $x$ such that the state of $x$ is flipped exactly at $p$ and $x$ at that moment is pivotal). Hence, we can sum over all sites $x$ in the small rhombus, and deduce that the conditional probability given $\Gamma^k[0, \sigma^k]$ that $t$ is very good is bounded from below by 
\begin{eqnarray*}
  \lefteqn {F^k ( \Gamma^k, t)+ (p_k-1/2) \times c'(\eta N_k)^2 \times P_{1/2} ( A^4 ( \eta N_k))} \\
  & \ge  & F^k ( \Gamma^k, t) + c'' \eta^2 \frac{P_{1/2} ( A^4 ( \eta N_k))}{P_{1/2} ( A^4 (N_k))} \times (p_k-1/2) N_k^2 P_{1/2} ( A^4 (N_k))\\
  &\ge & F^k ( \Gamma^k, t) + \frac { c''' \eta^2 }{P_{1/2} ( A^4 ( \eta N_k,N_k))}\\
  & \ge & F^k ( \Gamma^k, t) + \eta^{3/4 + \beta}
\end{eqnarray*}
(we have used the quasi-multiplicativity to estimate the ratio on the second line) for all small $\eta$ and uniformly with respect to $k$ and to the good configurations $\Gamma^k [0, \sigma^k]$. 
Letting $k \to \infty$,  we conclude that for any good triangle $t$, the conditional probability given $\gamma [0,\sigma]$ that $t$ is very good for $\gamma$ is at least 
$$ F( \gamma, t) + \eta^{3/4 + \beta}.$$

We now want to estimate the number of good triangles $t$. We divide the triangle $T$ into $M(\eta) \asymp \eta^{-2}$ small deterministic triangles. Using the same method as when studying the dimension of $\gamma$ (and the arm-separation ideas that allow to estimate the probability that a triangle is good i.e. to say that the probability for a triangle $t$ to be good is comparable to the probability to be hit), we see that with probability going to $1$ as $\eta \to 0$, the number of good triangles is at least $\eta^{-7/4 + \beta}$. 

We now count how many triangles are very good among the first $\eta^{-7/4+\beta}$ good ones, and we subtract the sum of the $F(t,\gamma)$ for these good triangles. If $\gamma$ was an SLE(6), the obtained quantity $Z(\eta)$ would satisfy $E( Z )  = 0$ and $\var ( Z) \le \eta^{-7/4 + \beta}$.
But, in our case, we see that
$$E (Z) \ge  \eta^{-7/4 + \beta} \times \eta^{3/4 + \beta } \ge \eta^{-1 + 2 \beta}$$
and $\var (Z) \le 4 \eta^{-7/4 + \beta}$.

If the sequence $(\eta_l)$ has been chosen to converge sufficiently fast to zero, we see that almost surely (using Borel-Cantelli), for all sufficiently large $l$, 
$$Z(\eta_l) \ge \eta_l^{-1+ 3 \beta},$$ whereas if $\gamma$ had been an SLE(6), we would almost surely have had 
$$Z(\eta_l) \le \eta_l^{-7/8 - \beta}$$ 
for infinitely many $l$'s.

We thus found an event that holds with probability one for $\gamma$, but with probability zero for SLE(6).
\qed
    
\section {Informal discussion}

Seemingly, the fact that the laws of near-critical interfaces are singular with respect to the law of SLE(6)
 surprises some theoretical physicists who work on this topic.
Recall that one important aspect of the SLE approach to critical systems was
precisely to show that critical conformally invariant models in the same
``universality class'' give rise to exactly the same curves in the scaling
limit. For near-critical models that are not strictly conformally invariant
such as near-critical percolation here, this strong ``universality'' can
fail to be true. 

Suppose that a random curve $\gamma$ is the scaling limit of a near-critical percolation interface. 
The tightness arguments, together with the estimates for existence of multiple arms indicate that the curve $\gamma$ is almost surely a Loewner chain, and that it can be defined via its random driving function $(w(t), t \ge 0)$. 
It is of course a natural question to ask what this driving function could be (recall that in the case of SLE(6), $w(t/6)$ is a standard Brownian motion). Note that after each (fixed or stopping) time, we expect the curve $\gamma$ to turn ``more right'' as what a SLE(6) would do in the same situation. This leads naturally to the conjecture that (if one extends properly the probability space) one can couple $(w(t), t \ge 0)$ with a standard Brownian motion $(\beta(t), t \ge 0)$ in such a way that $w(t) - \beta(6t)$ is an increasing continuous process adapted to the filtration generated by the processes $(\beta(6t))$ and $(w(t))$. 

One may think that this contradicts the fact that the law of $(w(t/6), t \ge 0)$ is singular with respect to that of $\beta$ (since the curves are generated by the driving functions, this is equivalent to say that the law of $\gamma$ is singular with respect to that of SLE(6)) because the law of Brownian motion with drift is known to be absolutely continuous with respect to that of Brownian motion. But this last fact is only valid under some regularity properties for the drift, and in the present case, one can expect this drift to be quite complicated (its derivative measure might be supported on a set of exceptional times i.e. a fractal-type set). 
%It seems possible using scaling argument what sort of implicit equation the driving function should satisfy, but it is not very enlightening -- at %least to us.   

We have seen that important properties (such as the dimension of the curve)
are shared by the critical interfaces and the near-critical interfaces in the
scaling limit, so that the technology based on conformal invariance of the
critical model still provides the correct description of the near-critical
interfaces in terms of exponents (this was used in {\cite{N1}} to describe the
``percolation front''). In fact, the near-critical percolation technology allows to prove the following result (this is a direct consequence of the quasi-multiplicativity property and the fact that arm probabilities are comparable uniformly to the arm probabilities at $p=1/2$):
for any $k$, there exists a constant $c_k$, such that for any $z_1, \cdots , z_k$ in $T$, and for any sufficiently small $\epsilon$, 
the ratio between the probability that $\gamma$ 
visits the $\epsilon$-neighborhood of these $k$ points 
and the probability of the same event for SLE(6) stays in $(c_k, 1/c_k)$ for all small $\epsilon$ and 
uniformly over the choices of $z_1, \ldots, z_k$. Loosely speaking, the ``finite-dimensional marginals'' of the law of  $\gamma$ 
are uniformly comparable to those of SLE(6) (even when the $z_i$'s are very close to each other). 
Note that this is the type of functions that is computable via the Conformal Field Theory technology. But this result does not 
contradict the fact that the law of $\gamma$ is singular with respect to that of SLE(6). Note for instance that the constant $c_k$ does depend on the 
number of considered points. Our proof used the fact that we did not zoom on the behavior of the curve at just one point, but that we detected the asymmetry by looking at the mean of the behavior of the curve near to more and more points.

\medbreak
\noindent
{\bf Acknowledgment.} We thank Christophe Garban for stimulating discussions.

\bigbreak
Laboratoire de Math\'ematiques, B\^at. 425, Universit\'e Paris-Sud 11, 91405 Orsay Cedex, France
\medbreak
D\'epartement de Math\'ematiques et Applications, \'Ecole Normale Sup\'erieure, 45 rue d'Ulm, 75230 Paris Cedex 05, France

\end{document}